\documentclass[12pt]{article}  

\usepackage{amsmath}
\usepackage{amssymb}
\usepackage{cite}
\usepackage{paralist}
\usepackage[usenames]{color}
\usepackage{graphicx}          
\usepackage[dvips]{epsfig}    

\title{Analysis of switching strategies for the optimization of periodic chemical reactions with controlled flow-rate}

\author{Peter Benner, Andreas Seidel-Morgenstern, Alexander Zuyev\footnote{Corresponding author is with the Max Planck Institute for Dynamics of Complex Technical Systems and is on leave from the Institute of Applied Mathematics and Mechanics, National Academy of Sciences of Ukraine.}}
\date{\small {Max Planck Institute for Dynamics of Complex Technical Systems}\\{Sandtorstra{\ss}e 1, 39106 Magdeburg, Germany}}

\begin{document}

\maketitle
\thispagestyle{empty}

\begin{abstract}
An isoperimetric optimal control problem with non-convex cost is considered for a class of nonlinear control systems with periodic boundary conditions. This problem arises in chemical engineering as the maximization of the product of non-isothermal reactions by consuming a fixed amount of input reactants. It follows from the Pontryagin maximum principle that the optimal controls are piecewise constant in the considered case. We focus on a parametrization of optimal controls in terms of switching times in order to estimate the cost under different switching strategies. We exploit the Chen--Fliess functional expansion of solutions to the considered nonlinear system with bang-bang controls to satisfy the boundary conditions and evaluate the cost analytically for small periods. In contrast to the previous results in this area, the system under consideration is not control-affine, and the integrand of the cost depends on the state.
This approach is applied to non-isothermal chemical reactions with simultaneous modulation of the input concentration and the volumetric flow-rate.
\end{abstract}

\section{Introduction}\label{sec_1}
Strategies for the dynamic optimization of chemical reaction models have been studied in the mathematical literature by using the Pontryagin maximum principle~\cite{BH1971,CES2017}, vibrational control technique~\cite{BBM1983}, frequency-domain methods~\cite{PS2013,SLZYW2018,SY1991}, center manifold theory~\cite{KDT2012}, flatness-based approach and extremum seeking~\cite{GDPH2007}, model predictive control methodology~\cite{ELC2017}, and other approaches.

A remarkable result in this area was formulated for a mathematical model of an isothermal reaction  the type ``$\nu_1 A_1 + \nu_2 A_2 \to$ Product''
with the power law rate $r=k C_1 ^{n_1} C_2^{n_2}$ in~\cite{Hoffmann}. Namely, it was shown that the conversion of $A_1$ and $A_2$ to the product cannot be improved by using time-varying controls if $0<n_1<1$, $0<n_2<1$, and $n_1+n_2\le 1$.
In the non-isothermal case, it turns out that it is possible to improve the performance of first-order reactions of the type ``$A \to$  Product'' by using sinusoidal periodic inputs~\cite{PS2013}.
For a realistic non-isothermal reaction of this type, it was shown that the optimal controls are bang-bang, and periodic switching strategies have been described  by applying the Pontryagin maximum principle in~\cite{CES2017}.
An analytic approach for computing the switching parameters of $\tau$-periodic controls has been developed in~\cite{AMM2019} for the case of small periods $\tau$.

Note that the above papers deal with reaction models with a constant flow-rate, while the periodic flow-rate modulation is shown to be an important ingredient for improving the reaction performance~\cite{NSP2014}.
The corresponding isoperimetric optimal control problem is rigorously formulated in~\cite{ZSB2019} for a non-isothermal mathematical model with two independent inputs: the inlet concentration and the flow-rate.
As in the case of constant flow-rate, it is shown in~\cite{ZSB2019} that the optimal controls are piecewise constant, and their switching times are defined in terms of zeros of certain auxiliary functions.
However, the structure of switching controllers has not been analyzed so far.
This paper aims at developing an efficient approach for computing periodic bang-bang controls and evaluating the cost for the isoperimetric optimal control problem introduced in~\cite{ZSB2019}.

\section{Optimization problem}

Consider a nonlinear control system describing non-isothermal chemical reactions of the type ``$A \to$  Product'' and order $\bar n$~\cite{NSP2014,ZSB2019}:
\begin{equation}
\dot x= f_0(x) + v_1 v_2 g_1(x) + v_2 g_2(x),\quad x=(x_1,x_2)^T\in{\mathbb R}^2,
 \label{CSTR_0}
\end{equation}
where $x_1$ is the dimensionless concentration of $A$ in the reactor, $x_2$ is the dimensionless temperature,
\begin{equation}
\begin{aligned}
f_0(x) &= \begin{pmatrix} -k_1(1+x_1)^{\bar n}\exp\{-\frac{\gamma}{x_2+1}\} \\ \delta -St(1+x_2) -k_2(1+x_1)^{\bar n}\exp\{-\frac{\gamma}{x_2+1}\} \end{pmatrix},\\
g_1(x) &=\begin{pmatrix}1+k_1 \exp\{-\gamma\} \\ 0\end{pmatrix},\; g_2(x)=\begin{pmatrix}-1-x_1 \\ k_2 \exp\{-\gamma\} + St - \delta - x_2\end{pmatrix},
\end{aligned}
\label{vector_fields}
\end{equation}
and $k_1$, $k_2$, $St$, $\gamma$, and $\delta$ are physical parameters (cf.~\cite{NSP2014}).
The dimensionless control variables $v_1\in [v_1^{min},v_1^{max}]$ and $v_2\in [v_2^{min},v_2^{max}]$ correspond to the inlet concentration of $A$ and the flow-rate, respectively.
We assume that $0<v_i^{min}\le 1$ and $v_i^{max}\ge 1$ for $i=1,2$.
Then it is easy to see that $x_1=x_2=0$ is an equilibrium of system~\eqref{CSTR_0} that corresponds to a steady-state operation of the considered chemical reactor with $v_1=v_2=1$.

System~\eqref{CSTR_0} can be transformed to the control-affine form with respect to the inputs  $u_1=v_1 v_2$ and $u_2=v_2$ as follows~\cite{ZSB2019}:
\begin{equation}
\dot x = f_0(x) + u_1 g_1(x) + u_2 g_2(x),\quad x\in {\mathbb R}^2,\; u=(u_1,u_2)^T\in U = {\rm Conv}\, U_b,
 \label{CSTR_affine}
\end{equation}
where
$$
U_b = \left\{\begin{pmatrix}u_1^{min} \\ u_2^{min}\end{pmatrix},\begin{pmatrix}u_1^{max} \\ u_2^{max}\end{pmatrix}, \begin{pmatrix}u_1^{-} \\ u_2^{max}\end{pmatrix}, \begin{pmatrix}u_1^{+} \\ u_2^{min}\end{pmatrix} \right\},
$$
$$
u_1^{min} = v_1^{min} v_2^{min},\; u_1^- = v_1^{min} v_2^{max},\; u_1^+ = v_1^{max} v_2^{min},\; u_1^{max} = v_1^{max} v_2^{max}.
$$

As maximizing the conversion of $A$ to the product over a given time period $t\in [0,\tau]$ can be treated in the sense of minimizing the remaning mass of $A$ in the outgoing stream,
our goal is to minimize the cost
\begin{equation}
 J =\frac{1}{\tau} \int_0^\tau \Bigl(x_1(t)+1\Bigr)u_2(t)dt.
\label{Jcost}
\end{equation}
We also assume that the consumption of $A$ over the period is fixed as $\frac{1}{\tau}\int_0^\tau u_1(t)dt=\bar u_1$, which yields the following isoperimetric optimal control problem.

{\bf Problem~2.1.}~\cite{ZSB2019}
{\em
Given $\tau>0$, $\bar u_1\in {\mathbb R}$, and $x^0\in {\mathbb R}^2$, the goal is to find an admissible control $\hat u\in L^\infty\Bigl([0,\tau];U\Bigr)$ that minimizes the cost $J$ along the trajectories of~\eqref{CSTR_affine}
corresponding to the admissible controls $u\in L^\infty\Bigl([0,\tau];U\Bigr)$ such that
\begin{equation}
\frac{1}{\tau}\int_0^\tau u_1(t)dt=\bar u_1\;\; \text{and}\;\; x(0)=x(\tau)=x^0.
\label{constraints}
\end{equation}
}

If $\hat u(t)$ $(0\le t\le \tau)$ is an optimal control for Problem~2.1, then it follows from the results of~\cite{ZSB2019} that $\hat u(t)\in U_b$ almost everywhere on $[0,\tau]$,
and the switching times of $\hat u(t)$ are related to zeros of the following functions:
$I_1(t)$ $I_2(t)$, $\frac{u_1^- - u_1^{min}}{u_2^{max} - u_2^{min}}I_1(t) + I_2(t)$, $\frac{ u_1^{max} - u_1^+}{u_2^{max} - u_2^{min}}I_1(t) + I_2(t)$,
where $I_1(t)$ and $I_2(t)$ are defined by solutions of the associated Hamiltonian system.
It should be noted that $I_1(t)$ and $I_2(t)$ are parameterized by initial values of the adjoint variables.
In this paper, we will not use any information on the behavior of adjoint variables and define the switching parameters directly from~\eqref{constraints}.
Then the cost~\eqref{Jcost} will be approximated analytically to estimate the performance improvement for the considered class of bang-bang controllers.

\section{Computation of the switching controls}

Assuming that a bang-bang control $\hat u(t)\in U_b$ $(0\le t\le \tau)$ has a finite number of switchings, we enumerate the switching times
\begin{equation}
0=t_0<t_1<...<t_N = \tau\;\; \text{with some}\;\; N\in \mathbb N
\label{tj}
\end{equation}
and denote
\begin{equation}
u^j = \hat u(t)\in U_b\;\; \text{for}\;\; t\in S_j = (t_{j-1},t_j),\; j=1,2,...,N.
\label{uj}
\end{equation}

Our goal is to analyse the cost $J$ on the trajectories of system~\eqref{CSTR_affine} with piecewise-constant controls of the form~\eqref{uj}
depending on the parameters $(t_1, ..., t_N)$ and $(u^1, ..., u^N)$.

A straightforward computation of $\int_0^\tau \hat u_1(t)dt$ for the piecewise-constant control~\eqref{uj} shows that the isoperimetric constraint in~\eqref{constraints} is equivalent to
\begin{equation}
\sum_{j=1}^N \alpha_j u^j_1 = \bar u_1\;\;\text{with}\;\; \alpha_j = \frac{t_j-t_{j-1}}{\tau}>0.
\label{iso_alpha}
\end{equation}
In order to satisfy the periodic boundary condition $x(0)=x(\tau)$ and estimate the cost~\eqref{Jcost} analytically for small $\tau$,
we exploit the Chen--Fliess expansion of solutions to system~\eqref{CSTR_affine} with the initial value $x(0)=x^0$ and control $u=\hat u (t)$ (see, e.g.,~\cite{AMM2019}):
\begin{equation}
\scriptsize
x = x^0 + \sum_{i=0}^2   g_i(x^0) V_i(t) + \sum_{i,j=0}^2\bigl(L_{g_j} g_i \bigr)(x^0) V_{ij}(t)  + \sum_{i,j,l=0}^2  \bigl( L_{g_l} L_{g_j} g_i \bigr) (x^0)V_{ijl}(t) + O(t^4),
\label{Fliess3}
\end{equation}
where we  assume that $g_0(x)=f_0(x)$, $L_{g_i}g_j(x)=\frac{\partial g_j(x)}{\partial x}g_i(x)$ is the directional derivative of $g_j(x)$ along $g_i(x)$, and
$$
\begin{aligned}
V_i(t)&=\int_0^t u_i(s )ds,\; u_0(t)\equiv 1,\; V_{ij}(t)=\int_0^t \int_0^s u_i(s)u_j(p)dp\,ds,\\ V_{ijl}(t)&=\int_0^t \int_0^s \int_0^p u_i(s)u_j(p)u_l(r)dr\, dp\,ds,\;\; t\in [0,\tau].
\end{aligned}
\label{Vall}
$$
The remainder of formula~\eqref{Fliess3} is of order $O(t^4)$ for small $t>0$ if the vector fields $g_j(x)$ are of class $C^3$ in a neighborhood of $x^0$.

As in~\cite{AMM2019}, we will restrict our analysis to the cases $N\le 4$,
motivated by the estimate of the number of switchings in isoperimetric problems proposed in~\cite{CES2017}.
The main analytical result of our study is summarized as follows.

{\bf Proposition~3.1.}
{\em
Let $\hat u(t)$, $t\in [0,\tau]$  be a bang-bang control represented by~\eqref{uj} with the parameters $0< t_1 \le t_2 \le t_3 \le t_4= \tau$ and $u^1,u^2,u^3,u^4\in U_b$,
and let $x(t)$, $t\in [0,\tau]$ be the corresponding solution of~\eqref{CSTR_affine} such that $x(0)=x^0\in {\mathbb R}^2$.
Then the isoperimetric constraint~\eqref{iso_alpha} is equivalent to
\begin{equation}
\sum_{j=2}^4 \alpha_j (u^j_1 - u_1^1) = \bar u_1 - u^1_1,\;\;\alpha_1 = 1-\alpha_2 - \alpha_3 -\alpha_4,
\label{iso_N4}
\end{equation}
and the periodic boundary condition $x(0)=x(\tau)$ reduces to
\begin{equation}
{\scriptsize
\begin{aligned}
&\sum_{j=1}^4 \alpha_j f_j  + \frac{\tau}{2}\bigl\{\alpha_1^2 L_{f_1}f_1+\alpha_2^2 L_{f_2}f_2 - \alpha_3^2 L_{f_3}f_3 - \alpha_4^2 L_{f_4}f_4+ 2 \alpha_1 \alpha_2 L_{f_1}f_2 -2 \alpha_3 \alpha_4 L_{f_4}f_3 \bigr\}
 \\
&  + \frac{\tau^2}{6}\bigl\{ \alpha_1^3 L_{f_1}^2 f_1 + \alpha_2^3 L_{f_2}^2 f_2 + \alpha_3^3 L_{f_3}^2 f_3  + \alpha_4^3 L_{f_4}^2 f_4 +3\alpha_1 \alpha_2 L_{f_1} (\alpha_1 L_{f_1}+\alpha_2 L_{f_2}) f_2\\
&+ 3\alpha_3 \alpha_4 L_{f_4} (\alpha_4 L_{f_4}+\alpha_3 L_{f_3}) f_3 \bigr\} = O(\tau^3),
\end{aligned}}
\label{periodic_N4}
\end{equation}
where $f_i(x)= f_0(x) + u^i_1 g_1(x) + u^i_2 g_2(x)$, $i=1,2,3,4$.
Moreover, the cost~\eqref{Jcost} evaluated for $x(t)$ admits the representation
$J= \bar u_2 + X_1$,
where
\begin{equation}
\bar u_2 = \frac{1}{\tau}\int_0^\tau \hat u_2(t)dt = u_2^1 + \sum_{j=2}^4 \alpha_j(u_2^j - u_2^1)
\label{u2_N4}
\end{equation}
and $X_1$ is the first component of the vector $X\in \mathbb R^2$:
\begin{equation}
\begin{aligned}
X=\frac{1}{\tau}&\int_0^\tau x(t) \hat u_2(t) \,dt = \bar u_2 x^0 + \frac{\tau}{2}\left(\alpha_1^2  u^1_2 f_1 - (1-\alpha_1)^2 u_2^2 f_2\right)\\
&+\frac{\tau^2}{6}\left( \alpha_1^3 u^1_2 L_{f_1}f_1+ (1-\alpha_1)^3 u^2_2 L_{f_2}f_2 \right)\\
&+\frac{\tau^3}{24}\left( \alpha_1^4 u^1_2 L_{f_1}L_{f_1}f_1 - (1-\alpha_1)^4 u^2_2 L_{f_2}L_{f_2}f_2 \right) + O(\tau^4).
\end{aligned}
\label{xu2}
\end{equation}
The vector fields $f_i(x)$ and their directional derivatives in~\eqref{periodic_N4},~\eqref{xu2} are evaluated at $x=x^0$.
}

The assertion of Proposition~3.1 is obtained from the Chen--Fliess expansion~\eqref{Fliess3} for the solution $x(t)$ of system~\eqref{CSTR_affine} with $u=\hat u(t)$.

Note that the cases with $N<4$ can be considered as particular cases of $N=4$ with some of the $\alpha_j$ being zero.
In particular, the case $N=2$ is treated by assuming $\alpha_3=\alpha_4=0$ in~\eqref{iso_alpha}. In this case, the equations~\eqref{iso_N4}, \eqref{periodic_N4}, and~\eqref{u2_N4}  are reduced, respectively, to
\begin{equation}
\alpha_1 = \frac{\bar u_1 - u^2_1}{u^1_1 - u_1^2}\in (0,1),\; \alpha_2= 1-\alpha_1\;\;\text{if}\;\; u^1_1 \neq u_1^2,
\label{iso_N2}
\end{equation}
{\scriptsize
\begin{equation}
\alpha_1 (f_1-f_2) + f_2 + \frac{\tau}{2}\left(\alpha_1^2 L_{f_1}f_1 - (1-\alpha_1)^2 L_{f_2}f_2\right) + \frac{\tau^2}{6}\left(\alpha_1^3 L_{f_1}^2 f_1 + (1-\alpha_1)^3 L_{f_2}^2 f_2\right) =O(\tau^3),
\label{periodic_N2}
\end{equation}}
and
\begin{equation}
\bar u_2 = \frac{1}{\tau}\int_0^\tau \hat u_2(t)dt =\alpha_1 u^1_2 + (1-\alpha_1)u^2_2.
\label{u2_N2}
\end{equation}

\section{Simulation results}
We take the following parameters for numerical simulations for the first-order $(\bar n=1)$ adiabatic reaction considered in~\cite{AMM2019}:
$$
\gamma = \frac{E_A}{R\bar T}=17.77,\, k_1 = k_0 {\bar C_A}^{\bar n-1} \frac{V}{\bar F}=5.819 \cdot 10^7,\,
k_2 = \frac{\Delta H_R k_0 {\bar C_A}^{\bar n} V}{\rho c_p \bar T \bar F}=-8.99\cdot 10^5,\,
\delta=St=0.
$$
The above dimensionless parameters are computed with the gas constant $$R= 8.3144598\, \frac{J}{K\cdot mol}$$ and the activation energy
$ E_A = 44.35\, \frac{kJ}{mol}$, the collision factor $k_0 = 1.4\cdot 10^5 \,s^{-1}$, the reaction heat $\Delta H_R = -55.5 \, \frac{kJ}{mol}$, and $\rho c_p = 4.186\, \frac{kJ}{K \cdot l}$ being the product of the density and the heat capacity.
This model corresponds to the chemical reaction $\rm (CH_3CO)_2 O + H_2 O \to 2\, CH_3 COOH$ in the CSTR of volume $V = 0.298\, l$ with the steady-state outlet concentration  ${\bar C_A}=0.3498\,\frac{mol}{l}$ and the steady-state temperature $\bar T = 300.17\, K$.
In contrast to the previous works~\cite{CES2017,AMM2019}, we consider the case of variable flow-rate in this paper.
Namely, we assume that the flow-rate and the inlet concentration can be controlled around their steady-state values $\bar F= 7.17\cdot 10^{-4}\,\frac{l}{s}$ and $\bar C_{Ai} = 0.74\, \frac{mol}{l} $, respectively, within the range of $85\%$, i.e.
$v_i^{min} = 0.15$, $v_i^{max} = 1.85$, $i=1,2$.
This choice of control constraints corresponds to the following components of the points in~$U_b$:
\begin{equation}
u_1^{min} = 0.0225,\; u_1^{max} = 3.4225,\; u_1^{+} = u_1^{-} = 0.2775,\;
u_2^{min} = 0.15,\; u_2^{max} = 1.85.
\label{control_constraints}
\end{equation}

In the sequel, we impose the isoperimetric constraint~\eqref{constraints} with $\bar u_1=1$.
The constraint $\bar u_1=1$ is satisfied, in particular, by the constant controls $u_1=u_2=1$ for system~\eqref{CSTR_affine} (or, equivalently, $v_1=v_2=1$ for system~\eqref{CSTR_0}).
As it was already mentioned, system~\eqref{CSTR_affine} admits the equilibrium $x_1=x_2=0$ with $u_1=u_2=1$,
and this equilibrium corresponds to the cost $\bar J =1$ in~\eqref{Jcost}.
In this section, we will compare the steady-state value $\bar J$ with the values of $J$ for the periodic trajectories corresponding to controls~\eqref{uj}.
As the goal of Problem~2.1 is to minimize the cost $J$, we will treat the periodic trajectories with $J<\bar J$ as improving the reactor performance in comparison with its steady-state operation.

\begin{figure}[h!]
{
\begin{minipage}[t]{.48\textwidth}
        \centering
        \includegraphics[width=1\linewidth,trim={2cm 14.8cm 8cm 2cm},clip]{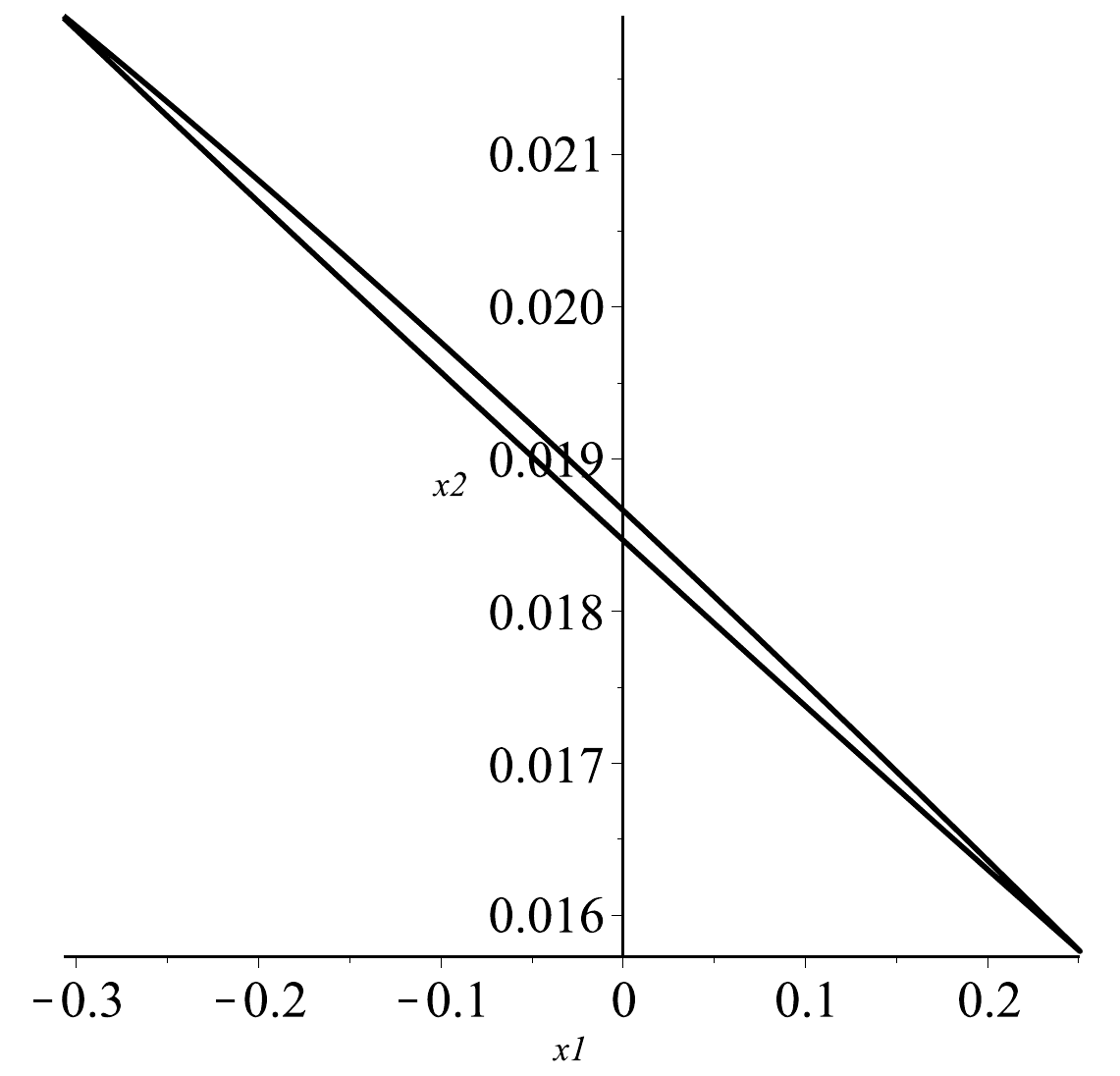}
        {\scriptsize a: Strategy~\eqref{C1}, $\tau=0.5$.}\\
        \vskip2ex
         \includegraphics[width=1\linewidth,trim={2cm 14.8cm 8cm 2cm},clip]{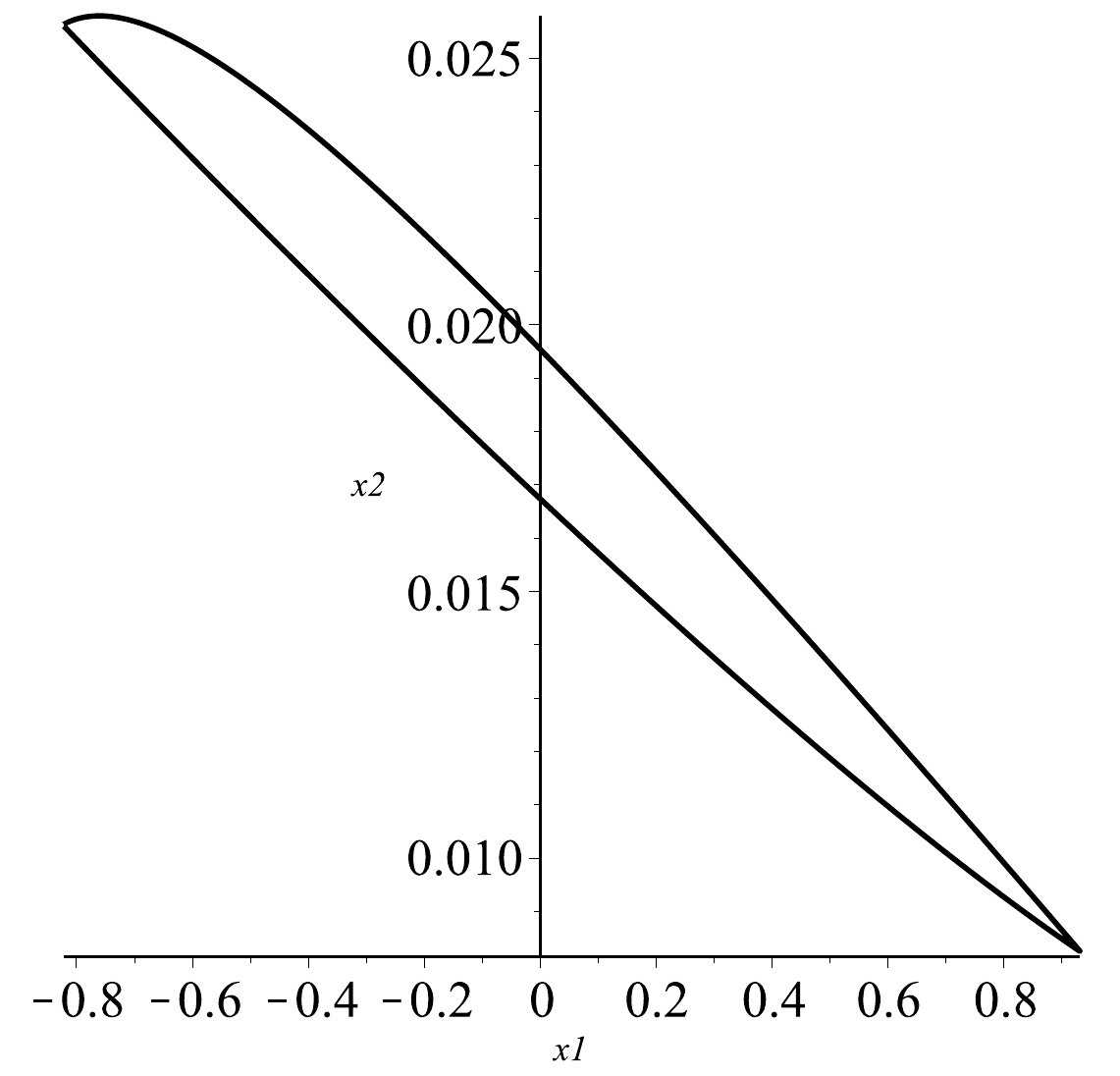}
        {\scriptsize c:  Strategy~\eqref{C1}, $\tau=2$.}
    \end{minipage}
    \hfill
    \begin{minipage}[t]{.48\textwidth}
        \centering
        \includegraphics[width=1\linewidth,trim={2cm 14.8cm 8cm 2cm},clip]{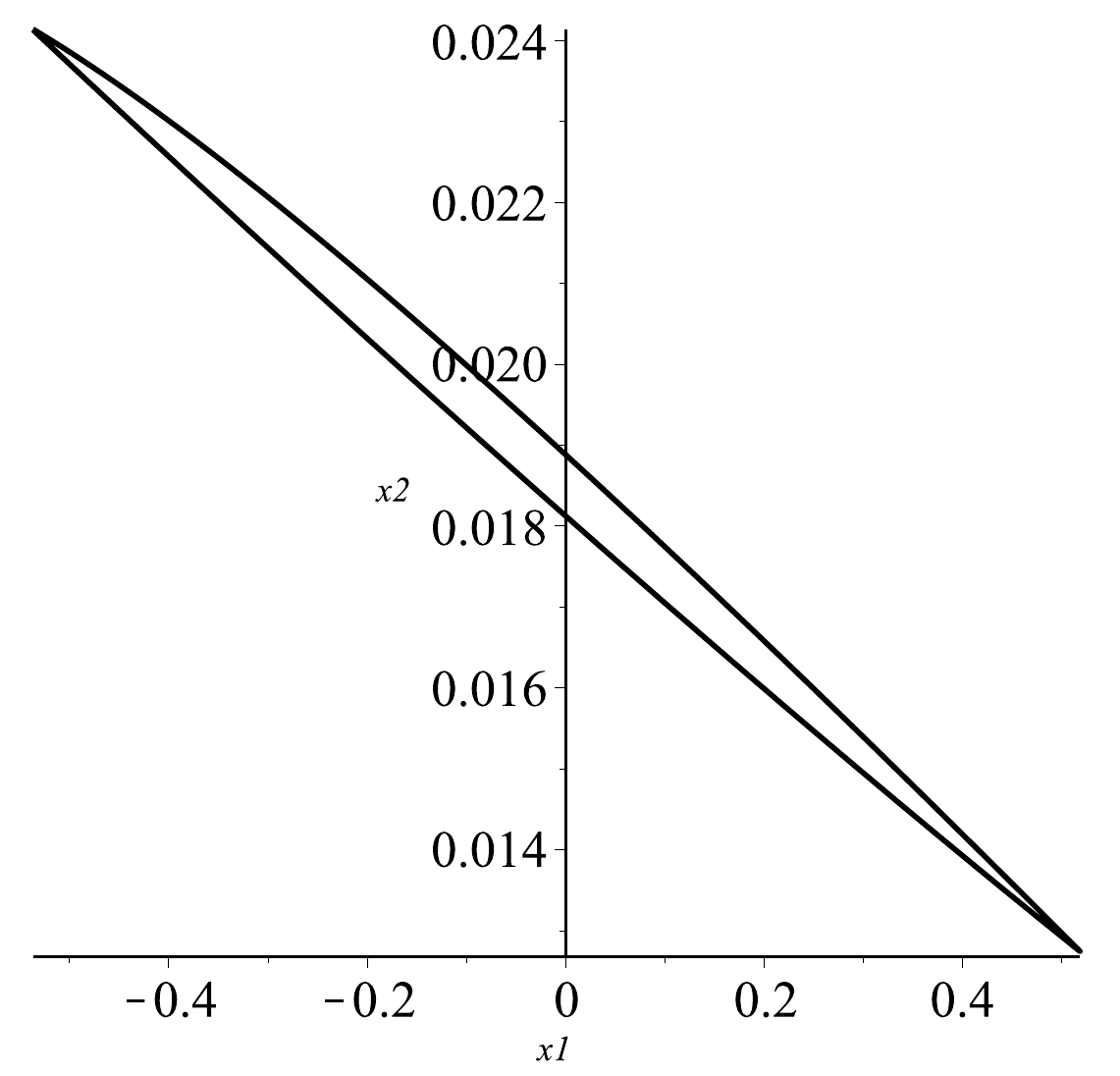}
        {\scriptsize b: Strategy~\eqref{C1}, $\tau=1$.}\\
        \vskip2ex
         \includegraphics[width=1\linewidth,trim={2cm 14.8cm 8cm 2cm},clip]{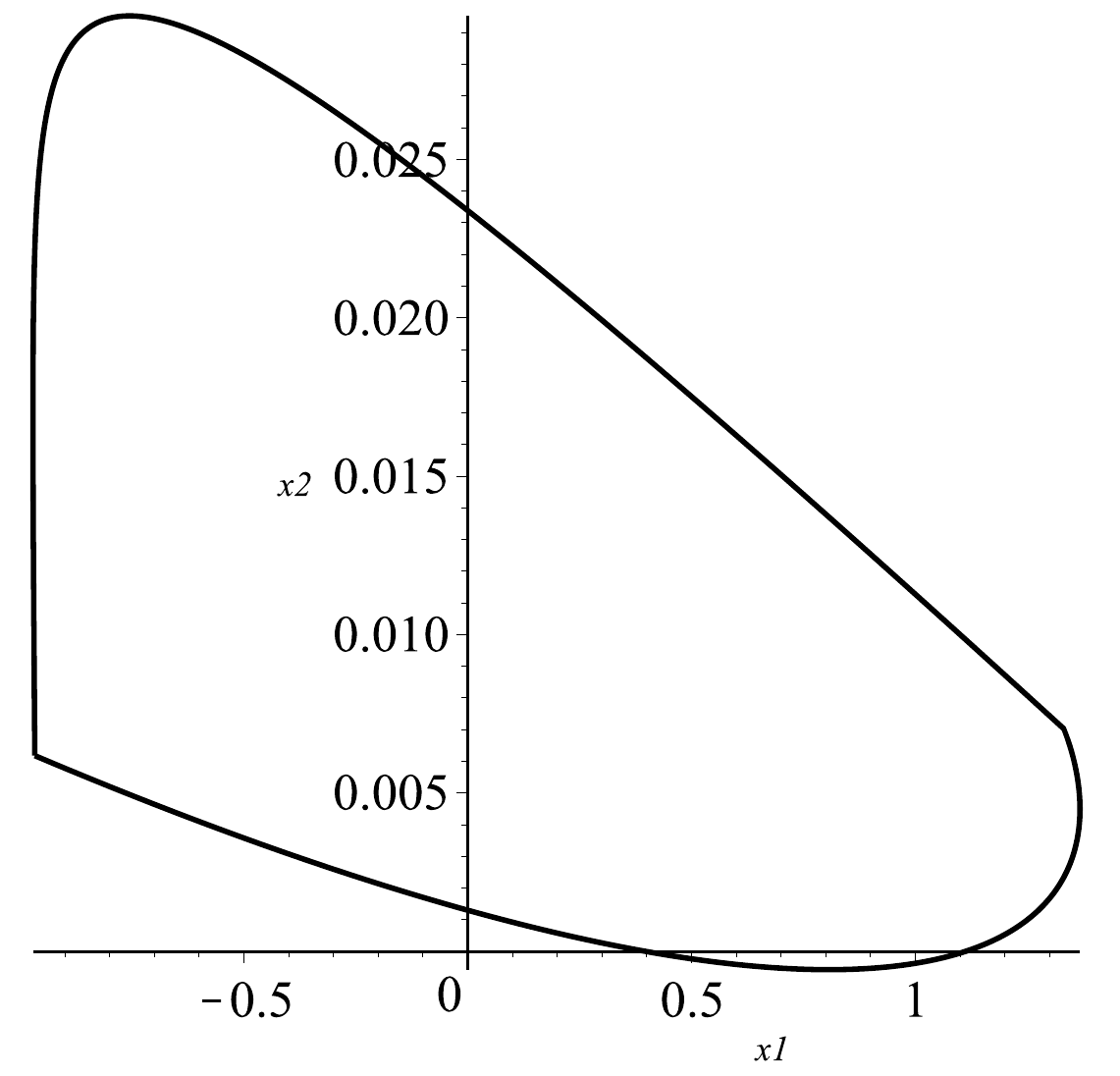}
        {\scriptsize d: Strategy~\eqref{C2}, $\tau=10$.}
    \end{minipage}
\caption{Periodic trajectories of system~\eqref{CSTR_affine} with $N=2$.}
\label{fig:1}}
\end{figure}

\begin{figure}[h!]
{
\begin{minipage}[t]{.48\textwidth}
        \centering
        \includegraphics[width=1\linewidth,trim={2cm 14.8cm 8cm 2cm},clip]{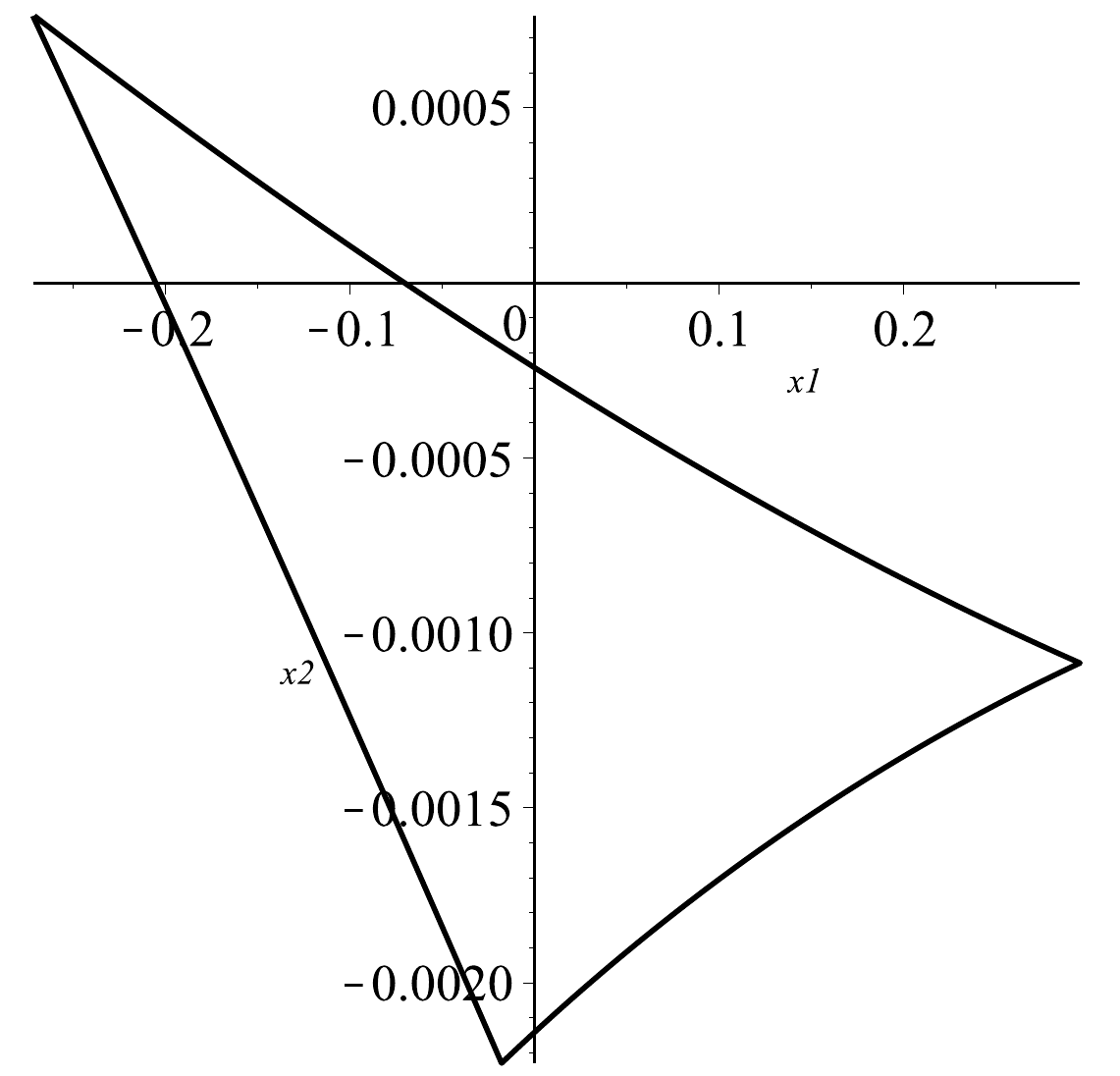}
        {\scriptsize a: Strategy~\eqref{C6}, $\tau=0.5$.}\\
        \vskip2ex
         \includegraphics[width=1\linewidth,trim={2cm 14.8cm 8cm 2cm},clip]{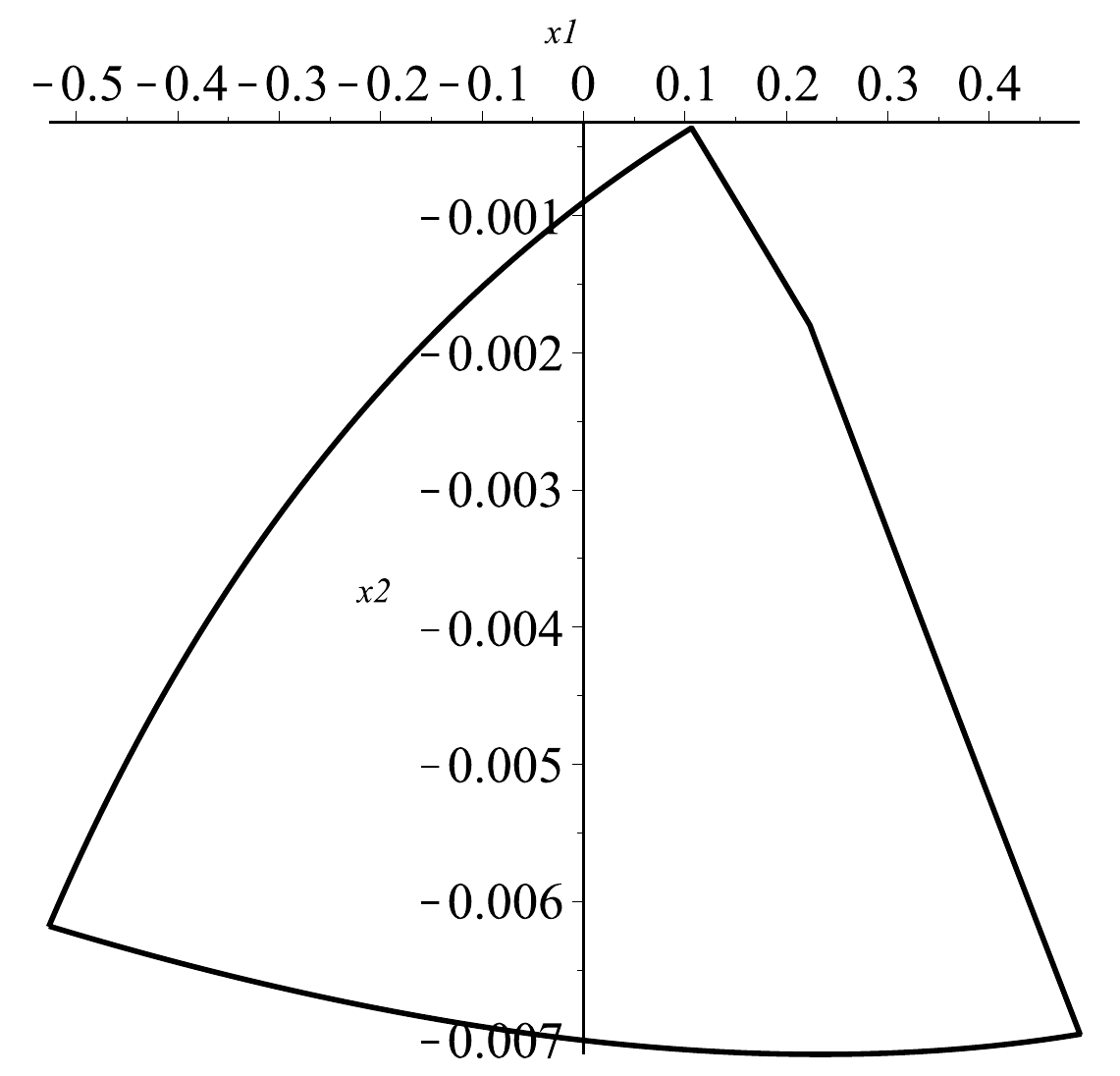}
        {\scriptsize c:  Strategy~\eqref{C7}, $\tau=1$.}
    \end{minipage}
    \hfill
    \begin{minipage}[t]{.48\textwidth}
        \centering
        \includegraphics[width=1\linewidth,trim={2cm 14.8cm 8cm 2cm},clip]{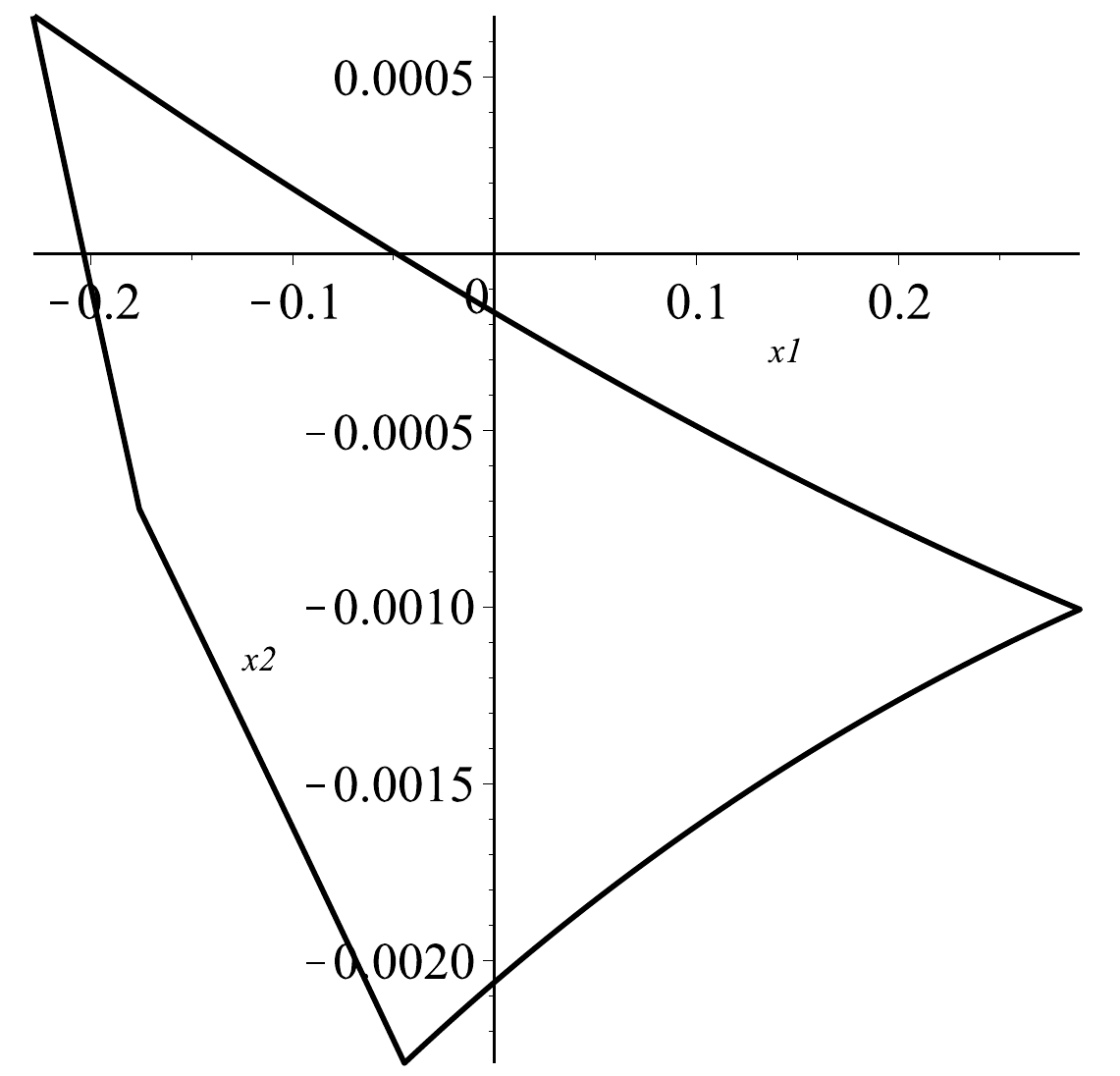}
        {\scriptsize b: Strategy~\eqref{C8}, $\tau=0.5$.}\\
        \vskip2ex
        \includegraphics[width=1\linewidth,trim={2cm 14.8cm 8cm 2cm},clip]{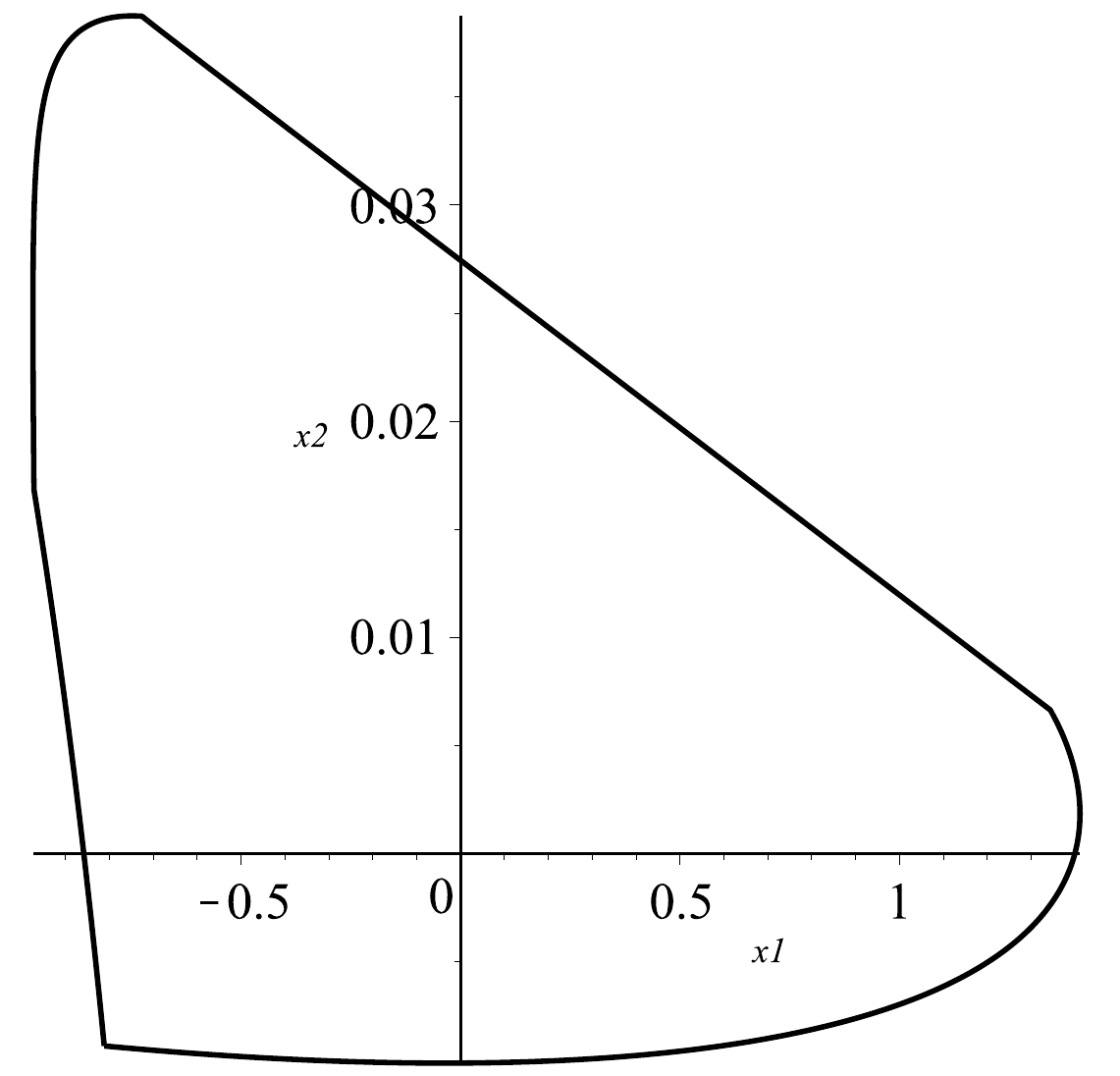}
        {\scriptsize d: Strategy~\eqref{C7}, $\tau=10$.}
    \end{minipage}
\caption{Periodic trajectories of system~\eqref{CSTR_affine} with $N=3$ and $N=4$.}
\label{fig:2}}
\end{figure}

The results of numerical simulations with controls of the form~\eqref{uj} are summarized in Table~1 and Figs.~\ref{fig:1}--\ref{fig:2} for the following switching strategies:
\begin{align}
&N=2,\; u^1 = \begin{pmatrix}
u_1^{max} \\ u_2^{max}
\end{pmatrix},\; u^2= \begin{pmatrix}
u_1^{min} \\ u_2^{min}
\end{pmatrix},\label{C1}\\
&N=2,\;u^1 = \begin{pmatrix}
u_1^{max} \\ u_2^{max}
\end{pmatrix},\; u^2 = \begin{pmatrix}
u_1^{+} \\ u_2^{min}
\end{pmatrix},\label{C2}\\
&N=3,\;u^1 = \begin{pmatrix}
u_1^{max} \\ u_2^{max}
\end{pmatrix},\; u^2 = \begin{pmatrix}
u_1^{min} \\ u_2^{min}
\end{pmatrix},
\; u^3 = \begin{pmatrix}
u_1^{-} \\ u_2^{max}
\end{pmatrix},
\label{C3}\\
&N=3,\;u^1 = \begin{pmatrix}
u_1^{max} \\ u_2^{max}
\end{pmatrix},\; u^2 = \begin{pmatrix}
u_1^{min} \\ u_2^{min}
\end{pmatrix},
\; u^3 = \begin{pmatrix}
u_1^{+} \\ u_2^{min}
\end{pmatrix},
\label{C4}\\
&N=3,\;u^1 = \begin{pmatrix}
u_1^{max} \\ u_2^{max}
\end{pmatrix},\; u^2 = \begin{pmatrix}
u_1^{+} \\ u_2^{min}
\end{pmatrix},
\; u^3 = \begin{pmatrix}
u_1^{-} \\ u_2^{max}
\end{pmatrix},
\label{C5}\\
&N=3,\;u^1 = \begin{pmatrix}
u_1^{max} \\ u_2^{max}
\end{pmatrix},\; u^2 = \begin{pmatrix}
u_1^{-} \\ u_2^{max}
\end{pmatrix},
\; u^3 = \begin{pmatrix}
u_1^{+} \\ u_2^{min}
\end{pmatrix},
\label{C6}\\
&N=4,\;
u^1 = \begin{pmatrix}
u_1^{max} \\ u_2^{max}
\end{pmatrix},\; u^2= \begin{pmatrix}
u_1^{+} \\ u_2^{min}
\end{pmatrix},\; u^3= \begin{pmatrix}
u_1^{min} \\ u_2^{min}
\end{pmatrix}, \; u^4=\begin{pmatrix}
u_1^{-} \\ u_2^{max}
\end{pmatrix},
\label{C7}\\
&N=4,\;
u^1 = \begin{pmatrix}
u_1^{max} \\ u_2^{max}
\end{pmatrix}, \; u^2=\begin{pmatrix}
u_1^{-} \\ u_2^{max}
\end{pmatrix}, \; u^3= \begin{pmatrix}
u_1^{min} \\ u_2^{min}
\end{pmatrix}, \; u^4=\begin{pmatrix}
u_1^{+} \\ u_2^{min}
\end{pmatrix}.
\label{C8}
\end{align}
Note that we only keep the switching strategies compatible with the constraint $\bar u_1=1$ in formulas~\eqref{C1}--\eqref{C8}, given the numerical values of controls in~\eqref{control_constraints}.
These formulas also allow the analysis of strategies obtained by cyclic permutations of $(u^1,u^2,u^3,u^4)$ because of the periodic nature of the considered control problem.
In Table~I, the switching parameters $\alpha_j=\frac{t_j-t_{j-1}}{\tau}$ are chosen according to the initial value $x^0$ of system~\eqref{CSTR_affine} by solving the algebraic equations~\eqref{iso_N4}, \eqref{periodic_N4} in Proposition~3.1.

\begin{table}[h!]
\begin{center}
\begin{tabular}{|c|c|c|c|c|}
\hline
{\bf Control} & {\bf Parameters} & {\bf Initial data} & {\bf Cost} \\ 
{\bf strategy} & {\bf $\alpha_j=(t_j-t_{j-1})/\tau$} & {\bf ${x^{0T}}$} & {\bf $J$} \\ \hline
\eqref{C1} & $\alpha_1=0.2875$, $\alpha_2=0.7125$ & $(-0.307, 0.0219)$ & {0.6293} \\ \hline
\eqref{C2} & $\alpha_1=0.2297$, $\alpha_2=0.7703$ & $(-0.3259, 0.0325)$ & {0.4883} \\ \hline
\eqref{C3} & $\alpha_1=0.2365$, $\alpha_2=0.0833$, $\alpha_3=0.6802$ & $(-0.2413, 0.017)$ & 0.653 \\ \hline
\eqref{C4} & $\alpha_1=0.2703$, $\alpha_2=0.5$, $\alpha_3=0.2297$ & $(-0.198, 0.00078)$ & 1.055 \\ \hline
\eqref{C5} & $\alpha_1=0.2297$, $\alpha_2=0.0833$, $\alpha_3=0.6870$ & $(-0.3305, 0.0312)$ & 0.502 \\ \hline
\eqref{C5} & $\alpha_1=0.2297$, $\alpha_2=0.1667$, $\alpha_3=0.6036$ & $(-0.3326, 0.0299)$ &  0.5169 \\ \hline
\eqref{C5} & $\alpha_1=0.2297$, $\alpha_2=0.25$, $\alpha_3=0.5203$ & $(-0.332, 0.0287)$ & 0.5326 \\ \hline
\eqref{C5} & $\alpha_1=0.2297$, $\alpha_2=0.3333$, $\alpha_3=0.4370$ & $(-0.3306, 0.0273)$ & 0.5488 \\ \hline
\eqref{C5} & $\alpha_1=0.2297$, $\alpha_2=0.4167$, $\alpha_3=0.3536$ & $(-0.3269, 0.026)$ & 0.5659 \\ \hline
\eqref{C5} & $\alpha_1=0.2297$, $\alpha_2=0.5$, $\alpha_3=0.2703$ & $(-0.323, 0.0249)$ & 0.5828 \\ \hline
\eqref{C6} & $\alpha_1=0.2297$, $\alpha_2=0.5$, $\alpha_3=0.2703$  & $(-0.271, 0.00076)$ & 1.0591 \\ \hline
\eqref{C7} & $\alpha_1=0.264$, $\alpha_2=0.083$, $\alpha_3=0.417$, $\alpha_4=0.236$ & $(-0.329, -0.0056)$ & 1.1259 \\ \hline
\eqref{C7} & $\alpha_1=0.237$, $\alpha_2=0.417$, $\alpha_3=0.083$, $\alpha_4=0.263$  & $(-0.263, 0.0133)$ & 0.7179 \\ \hline
\eqref{C7} & $\alpha_1=\alpha_2=\alpha_3=\alpha_4=0.25$ & $(-0.266, 0.00066)$ & 0.9465 \\ \hline
\eqref{C8} & $\alpha_1=0.264$, $\alpha_2=0.083$, $\alpha_3=0.417$, $\alpha_4=0.236$  & $(-0.2077, 0.0007)$ & 1.057 \\ \hline
\eqref{C8} & $\alpha_1=0.237$, $\alpha_2=0.417$, $\alpha_3=0.083$, $\alpha_4=0.263$& $(-0.256, 0.0007)$ & 1.0604 \\ \hline
\eqref{C8} & $\alpha_1=\alpha_2=\alpha_3=\alpha_4=0.25$ & $(-0.228, 0.00067)$ & 1.0616 \\ \hline
\end{tabular}
\caption{Simulation results for system~\eqref{CSTR_affine} with controls~\eqref{uj}, $\tau=0.5$.}
\end{center}
\label{tab:01}
\end{table}

\section{Conclusions}
The presented simulation results confirm that the best performance improvement in the sense of the cost~\eqref{Jcost} is achieved by bang-bang controls of the form~\eqref{uj} in the case~\eqref{C2} (up to a permutation of $u^1$ and $u^2$).
Note that the periodic trajectories in Figs.~\ref{fig:1} and~\ref{fig:2} are obtained as numerical solutions of system~\eqref{CSTR_affine},~\eqref{uj}, and their orbital stability (or partial stability~\cite{Z2015}) remains to be verified in future work to justify the practical relevance of the proposed discontinuous control strategies.


\begin{thebibliography}{10}
\bibitem{BH1971}
{\sc Bailey, J., and Horn, F.}
\newblock Comparison between two sufficient conditions for improvement of an
  optimal steady-state process by periodic operation.
\newblock {\em Journal of Optimization Theory and Applications 7}, 5 (1971),
  378--384.

\bibitem{BBM1983}
{\sc Bellman, R., Bentsman, J., and Meerkov, S.}
\newblock Vibrational control of systems with Arrhenius dynamics.
\newblock {\em Journal of Mathematical Analysis and Applications 91}, 1 (1983),
  152--191.

\bibitem{AMM2019}
{\sc Benner, P., Seidel-Morgenstern, A., and Zuyev, A.}
\newblock Periodic switching strategies for an isoperimetric control problem
  with application to nonlinear chemical reactions.
\newblock {\em Applied Mathematical Modelling 69\/} (2019), 287--300.

\bibitem{ELC2017}
{\sc Ellis, M., Liu, J., and Christofides, P.}
\newblock {\em Economic Model Predictive Control: Theory, Formulations and
  Chemical Process Applications}.
\newblock Springer, London, 2017.

\bibitem{GDPH2007}
{\sc Guay, M., Dochain, D., Perrier, M., and Hudon, N.}
\newblock Flatness-based extremum-seeking control over periodic orbits.
\newblock {\em IEEE Transactions on Automatic Control 52}, 10 (2007),
  2005--2012.

\bibitem{Hoffmann}
{\sc Hoffmann, U., and Sch\"adlich, H.-K.}
\newblock The influence of reaction orders and of changes in the total number
  of moles on the conversion in a periodically operated CSTR.
\newblock {\em Chemical Engineering Science 41\/} (1986), 2733--2738.

\bibitem{KDT2012}
{\sc Kravaris, C., Dermitzakis, I., and Thompson, S.}
\newblock Higher-order corrections to the pi criterion using center manifold
  theory.
\newblock {\em European Journal of Control 18}, 1 (2012), 5--19.

\bibitem{NSP2014}
{\sc Nikoli\'c, D., Seidel-Morgenstern, A., and Petkovska, M.}
\newblock Nonlinear frequency response analysis of forced periodic operation of
  non-isothermal cstr using single input modulations. Part~I: Modulation of
  inlet concentration or flow-rate.
\newblock {\em Chemical Engineering Science 117\/} (2014), 71–84.

\bibitem{PS2013}
{\sc Petkovska, M., and Seidel-Morgenstern, A.}
\newblock Evaluation of periodic processes.
\newblock In {\em Periodic Operation of Chemical Reactors}, P.~Silveston and
  R.~Hudgins, Eds. Butterworth-Heinemann, 2013, p.~387–413.

\bibitem{SLZYW2018}
{\sc Shi, H., Lang, Z., Zhu, Y., Yuan, D., and Wang, W.}
\newblock Optimal design of the inlet temperature based periodic operation of
  non-isothermal CSTR using nonlinear output frequency response functions.
\newblock {\em IFAC-PapersOnLine 51}, 18 (2018), 620--625.

\bibitem{SY1991}
{\sc Sterman, L., and Ydstie, B.}
\newblock Periodic forcing of the cstr: An application of the generalized
  $\pi$-criterion.
\newblock {\em AIChE Journal 37}, 7 (1991), 986--996.

\bibitem{Z2015}
{\sc Zuyev, A.}
\newblock {\em Partial Stabilization and Control of Distributed Parameter
  Systems with Elastic Elements}.
\newblock Springer, Cham, 2015.

\bibitem{CES2017}
{\sc Zuyev, A., Seidel-Morgenstern, A., and Benner, P.}
\newblock An isoperimetric optimal control problem for a non-isothermal
  chemical reactor with periodic inputs.
\newblock {\em Chemical Engineering Science 161\/} (2017), 206--214.

\bibitem{ZSB2019}
{\sc Zuyev, A., Seidel-Morgenstern, A., and Benner, P.}
\newblock Optimal periodic control of nonlinear chemical reactions with a
  time-varying flow rate.
\newblock {\em PAMM - Proceedings in Applied Mathematics and Mechanics 19\/}
  (2019), in Press, DOI: 10.1002/pamm.201900160.
  \end{thebibliography}
\end{document}